\documentclass{article}

\usepackage{graphicx}
\usepackage{amssymb}
\usepackage{amsmath}
\usepackage{amsthm}
\usepackage{tikz}
\usepackage{hyperref}
\usepackage[hang,flushmargin]{footmisc}
\usepackage{colortbl}

\usetikzlibrary{positioning}

\theoremstyle{plain}
\newtheorem{definition}{Definition}
\newtheorem{theorem}{Theorem}

\newtheorem{conjecture}{Conjecture}

\newcommand\blfootnote[1]{%
  \begingroup
  \renewcommand\thefootnote{}\footnote{#1}%
  \addtocounter{footnote}{-1}%
  \endgroup
}

\def\irr#1{{\rm Irr}(#1)}

\begin{document}
\newpage

\title{A Divisibility Problem in the McKay Conjecture}
\author{Ignacio P. Navarro}
\maketitle
\begin{abstract}
The mathematical software \texttt{GAP} (Groups, Algorithms, Programming)
offers a powerful set of tools to investigate computationally group theory.
Using this software package we investigate a variation of a well-known 
problem in representation theory, the McKay Conjecture,
that asserts that there is a bijection between two sets
of complex-valued functions defined on two finite groups. In order to do so,
we use the extensive group database from \texttt{GAP} and Kuhn's Algorithm
in order to test our conjecture. In the important case of symmetric groups, 
some number theory might help to explain our findings.
\end{abstract}
\blfootnote{\emph{Imperial College London} \\ Email: \texttt{ignacio.navarro14@imperial.ac.uk}}

\section{Preliminaries}
One of the most important operations in a group $G$, besides the multiplication of elements in 
the group, is conjugation: if $g, x \in G$, then we can \emph{conjugate} $G$ by $x$, and this
is $g^x = x^{-1}gx$. In this case, $g$ and $g^x$ are \textbf{conjugate} in $G$. Conjugation
appears everywhere in Mathematics, perhaps the best example is whenever we have a linear
map $f$ in a vector space and we find the matrix of $f$ with respect to two bases.
The relation $g \equiv h$ if and only if $g$ and $h$ are $G$-conjugate defines an equivalence relation in $G$, which partitions it into \textbf{conjugacy classes}.
The conjugacy class of $g \in G$ is simply the set $\{ g^x \mid x \in G \}$ of all
conjugates of $g$. In an abelian group, of course $g^x = x^{-1}gx = gx^{-1}x = g$,
and the conjugacy classes of $G$ are simply the group elements. This is the less
interesting case.

\medskip

In a group, we conjugate not only elements but subgroups: if $H$ is a subgroup
of $G$ and $x \in G$, then
$$
H^x = \{ x^{-1}hx \mid h \in H \}
$$
is another subgroup of $G$. Furthermore, $H$ is \textbf{normal} in $G$, and we
write $H \trianglelefteq G$, if $H^x = H$ for all $x \in G$.

\medskip

Now, we consider the complex space of functions $f: G \rightarrow \mathbb{C}$ 
(a space of dimension $|G|$), and we are interested in the functions
that satisfy $$f(g^x)=f(g)$$ for all $g, x \in G$, i.e., the functions that
are constant on conjugacy classes. This subspace $\text{cf}(G)$ of the 
\textbf{class functions} of $G$ has dimension the number of conjugacy classes of $G$. This is easy to prove: the \emph{characteristic functions} $f_K$ that take the value 0 outside the conjugacy class $K$ and take the value 1 on the conjugacy class $K$ form a basis of $\text{cf}(G)$.
As a good example of the functions that we are interested, consider for instance
the group of invertible matrices $\text{GL}_n(\mathbb{C})$ and consider $f$ to be the trace function. Since conjugate matrices have the same trace, then $f$ is a class function.
    
\medskip

More generally, if $G$ is any finite group, then we are interested in the
\textbf{representations} of $G$, which are the
multiplicative functions $\mathcal{X}: G \rightarrow \text{GL}_n(\mathbb{C})$, and their trace functions.
These trace functions $f: G \rightarrow \mathbb{C}$ are called the \textbf{characters} of $G$.
These functions are crucial in finite group theory.
Observe that if $\chi$ is the trace of the representation $\mathcal{X}$, then
$\chi(1)$ is the trace of the identity function, and therefore $\chi(1) = n$. This
is called the \textbf{degree} of the character.

\medskip

A new way to construct representations is by using \emph{direct sums}: 
if $\mathcal{Y}$ and $\mathcal{Z}$ are representations of $G$, then we can form
a new representation
$$ 
\left( 
\begin{array}{cc}
  \mathcal{Y}  &     0       \\
       0       & \mathcal{Z} \\
\end{array} 
\right)
$$
whose character is the sum of the characters given by $\mathcal{Y}$ and $\mathcal{Z}$.

\medskip

\begin{definition}
A character $\chi$ of $G$ is irreducible
if it is not the sum of two characters.
\end{definition}

With the previous definition, we now come to the Fundamental Theorem of Character Theory.

\begin{theorem}
If $G$ is a finite group, then the set $\irr{G}$ of the irreducible characters of $G$ form
a basis of the space of class functions. In particular, 
$|\irr{G}|$ is the number of conjugacy classes of $G$.
\end{theorem}

\section{The McKay Conjecture}

In order to state the McKay conjecture, we need a bit more of group theory.
We fix a prime $p$ and a finite group $G$, and we consider
a \textbf{Sylow $p$-subgroup} of $G$, which is a subgroup $P \le G$ of order
the largest power of $p$ that divides $|G|$. These subgroups exist
(by Sylow theorems), and if $P$ and $Q$ are Sylow subgroups of $G$, then there is $x \in G$
such that
$$
P^x=\{ x^{-1}y x \mid y \in P \}=Q
$$

The subgroup of $G$ consisting of $\{x \in G \mid P^x=P \} = \textbf{N}_{G}(P)$,
is called the \textbf{normalizer} of $P$ in $G$.

\medskip

Essentially, the McKay conjecture states that an important number associated with $G$ can actually be computed in a much smaller group, $\textbf{N}_{G}(P)$. This was formulated in 1971 in by J. McKay and asserts the following:

\smallskip

\begin{conjecture}
If $G$ is a finite group, $p$ is a prime and $\text{Irr}_{p'}(G)$ is the set
of the irreducible characters of $G$ of degree not divisible by $p$ then
there is a bijection 
$$
f:{\rm Irr}_{p'}(G) \rightarrow {\rm Irr}_{p'}(\textbf{N}_{G}(P))
$$
In other words,
$$|{\rm Irr}_{p'}(G)|=|{\rm Irr}_{p'}(\textbf{N}_{G}(P))|\, .$$
\end{conjecture}

\medskip

Our aim in this note is to investigate using the mathematical software \texttt{GAP} when it is possible to choose bijections
$$f:{\rm Irr}_{p'}(G) \rightarrow {\rm Irr}_{p'}(\textbf{N}_{G}(P)) $$
such that 
$$f(\chi)(1)\, \text{ divides } \, \chi(1)$$
for all $\chi \in \irr G$.  This seems not to have been investigated until now.
It also presents an interesting computational problem that we shall solve using a bipartite
graph and the Hungarian Algorithm. Our findings allow us to propose the following conjecture. If $p$ is a prime, first recall that a finite group $G$ is \textbf{$p$-solvable} if it has normal subgroups
$$1=G_0 \triangleleft G_1 \triangleleft \ldots \triangleleft G_{n-1} \triangleleft G_n=G$$
such that $|G_i|/|G_{i-1}|$ has either order a power of $p$ or not divisible by $p$.
In particular, solvable groups are $p$-solvable for every prime $p$.     

\medskip

\begin{conjecture}\label{conjsolvable}
Let $G$ be a $p$-solvable group. Then there exists a bijection
$$f:{\rm Irr}_{p'}(G) \rightarrow {\rm Irr}_{p'}(\textbf{N}_{G}(P)) $$
such that 
$$f(\chi)(1)\, \text{ divides } \, \chi(1) \quad \forall \chi \in \irr G.$$
\end{conjecture}  

\section{Algorithm and \texttt{GAP} Database}
To test the conjecture, we will find the bijection explicitly for most of the groups (memory permitted) in \texttt{GAP}. 
Indeed we have two lists (and not sets, as we are looking for a bijection)
\begin{align*}
A &= [ \; \chi(1) \; | \;  \chi \in \text{Irr}_{p'}(G) \;] \\
B &= [ \; \psi(1) \; | \;  \psi \in \text{Irr}_{p'}(\textbf{N}_{G}(P)) \;]
\end{align*}

And we are trying to find a bijection $f:A \to B$ such that $f(a)$ divides $a$ for all $a \in A$.
This is equivalent to finding a permutation $B'$ of $B$ such that $B'[i]$ divides $A[i]$ for all $i$.
The straightforward approach to find the bijection would be to find all the permutations
of $B$ and see if one satisfies our divisibility condition. This however poses a computational problem for big groups.
For example, in
the symmetric group $S_{25}$ and $p=5$, there are 25 irreducible characters of
degree not divisible by 5. That means that
list $B$ has 25 elements, and then we are looking at $25! = 15511210043330985984000000$ permutations at worst.

\medskip

There is, however, a more elegant approach lowering the complexity from $O(n!)$ to $O(n^3)$.
Let $G = (A, B, E)$ be a bipartite graph. Connect each $a \in A$ with all the divisors
from $B$. Now see if the maximum matching in $G$ is a perfect matching, where the number of
edges is equal to the number of pairs $(a, b)$ that are connected. This combinatorial
optimization algorithm is called the Hungarian Algorithm, and although it is used to solve
the \textbf{assignment problem}, we use it in our advantage by setting throughout the weight function
as $W \equiv 1$. 

\begin{figure}[!h]
\centering
\begin{tikzpicture}[every node/.style={circle,draw}]

	\node[draw=none] at (-2,6) {$A$};
	\node[draw=none] at (7,6) {$B$};

    \node[minimum size=1.66cm, inner sep=0pt] at (0, 9) (a1) {\footnotesize$\chi_1(1) = 6$};
    \node[minimum size=1.66cm, inner sep=0pt] at (0, 7) (a2) {\footnotesize$\chi_2(1) = 10$};
    \node[minimum size=1.66cm, inner sep=0pt] at (0, 5) (a3) {\footnotesize$\chi_3(1) = 22$};
    \node[minimum size=1.66cm, inner sep=0pt] at (0, 3) (a4) {\footnotesize$\chi_4(1) = 26$};
    \node[minimum size=1.66cm, inner sep=0pt] at (5, 9) (b1) {\footnotesize$\psi_1(1) = 3$};
    \node[minimum size=1.66cm, inner sep=0pt] at (5, 7) (b2) {\footnotesize$\psi_2(1) = 5$};
    \node[minimum size=1.66cm, inner sep=0pt] at (5, 5) (b3) {\footnotesize$\psi_3(1) = 13$};
    \node[minimum size=1.66cm, inner sep=0pt] at (5, 3) (b4) {\footnotesize$\psi_4(1) = 2$};

    \draw[color={rgb:red,21;green,133;blue,21}] (a1) -- (b1); \draw(a1) -- (b4);
    \draw[color={rgb:red,21;green,133;blue,21}](a2) -- (b2); \draw(a2) -- (b4);
    \draw[color={rgb:red,21;green,133;blue,21}](a3) -- (b4);
    \draw[color={rgb:red,21;green,133;blue,21}](a4) -- (b3); \draw(a4) -- (b4);

\end{tikzpicture}
\caption{Example bipartite graph to find the bijection using Kuhn's Algorithm}
\label{graph}
\end{figure}
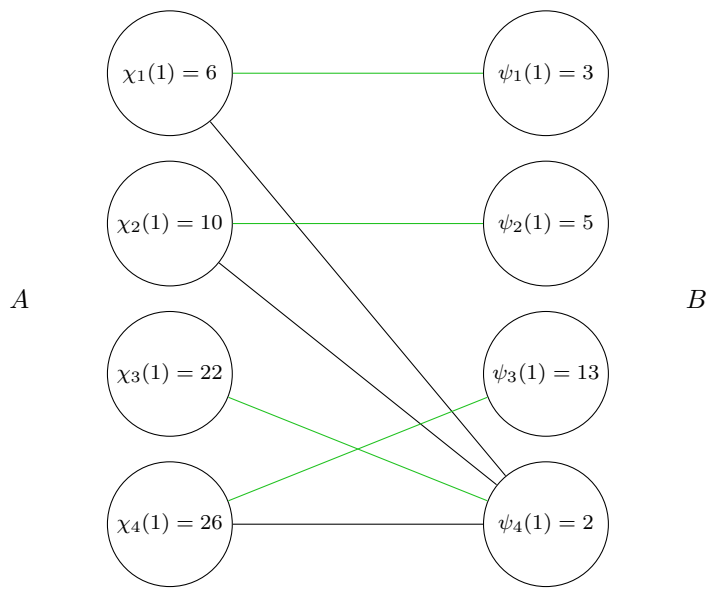

For example, in Figure~\ref{graph}, we see that there is in fact a bijection (labeled with green edges), and therefore
a permutation $B'$ on $B$ such that $B'[i]$ divides $A[i]$ for all $i$.

\medskip

The \texttt{GAP} package contains all groups of order up to 2000 (and even more for some cases) except the groups
of order 1024 (there are too many of those). Each group is labeled by $\texttt{SmallGroup}(x,y)$, where this means that $\texttt{SmallGroup}(x,y)$ is the $y$-th group of order $x$. We will use this 
in our investigation.

\section{Results}

\subsection{Results for solvable groups}

The results for solvable groups motivate Conjecture~\ref{conjsolvable}.

\subsection{Results for symmetric groups}

\begin{table}[]
\centering
\def\arraystretch{1.5}
\begin{tabular}{|c|c|c|c|c|c|c|c|c|}
\hline
       & $S_2$ & $S_3$ & $S_5$ & $S_7$ & $S_{11}$ & $S_{13}$ & $S_{17}$  & $S_{19}$               \\ \hline
$p=2$  & \cellcolor{blue!25}\checkmark & \checkmark & \checkmark & \checkmark  & \checkmark & \checkmark & \checkmark & \checkmark \\ \hline
$p=3$  & \cellcolor{green!25}\checkmark & \cellcolor{blue!25}\checkmark & \checkmark & \text{\sffamily{X}} & \checkmark  & \checkmark   & \checkmark  & \checkmark           \\ \hline
$p=5$  & \cellcolor{green!25}\checkmark & \cellcolor{green!25}\checkmark & \cellcolor{blue!25}\checkmark & \text{\sffamily{X}} & \text{\sffamily{X}} & \text{\sffamily{X}} & \text{\sffamily{X}} & \checkmark           \\ \hline
$p=7$  & \cellcolor{green!25}\checkmark & \cellcolor{green!25}\checkmark & \cellcolor{green!25}\checkmark & \cellcolor{blue!25}\checkmark   & \checkmark  & \text{\sffamily{X}} & \checkmark   & \checkmark    \\ \hline
$p=11$ & \cellcolor{green!25}\checkmark & \cellcolor{green!25}\checkmark & \cellcolor{green!25}\checkmark & \cellcolor{green!25}\checkmark   & \cellcolor{blue!25}\checkmark  & \checkmark   & \text{\sffamily{X}} & \text{\sffamily{X}} \\ \hline
$p=13$ & \cellcolor{green!25}\checkmark & \cellcolor{green!25}\checkmark & \cellcolor{green!25}\checkmark & \cellcolor{green!25}\checkmark   & \cellcolor{green!25}\checkmark  & \cellcolor{blue!25}\checkmark   & \checkmark    & \text{\sffamily{X}} \\ \hline
$p=17$ & \cellcolor{green!25}\checkmark & \cellcolor{green!25}\checkmark & \cellcolor{green!25}\checkmark & \cellcolor{green!25}\checkmark   & \cellcolor{green!25}\checkmark  & \cellcolor{green!25}\checkmark   & \cellcolor{blue!25}\checkmark    & \checkmark           \\ \hline
$p=19$ & \cellcolor{green!25}\checkmark & \cellcolor{green!25}\checkmark & \cellcolor{green!25}\checkmark & \cellcolor{green!25}\checkmark   & \cellcolor{green!25}\checkmark  & \cellcolor{green!25}\checkmark   & \cellcolor{green!25}\checkmark    & \cellcolor{blue!25}\checkmark           \\ \hline
\end{tabular}
\caption{Whether a bijection was found for different primes and symmetric groups. Cells in green are trivial cases, while cells in blue are the interesting case.}
\label{symmetrictable}
\end{table}

When investigating symmetric groups, there was no clear pattern to start off. However, 
after fixing the first few primes and seeing the relationship, we can start to see a pattern.

\medskip 

Indeed, fix an ordered set of primes $P = \{ p_1, p_2, ..., p_n \}$ where $p_i < p_j$ if $i < j$.
Set $S = \{ S_{p_1}, S_{p_2}, ..., S_{p_n} \}$, where $S_{p_i}$ is the Symmetric group on $p_i$ elements. Finally, set $M$ to be a $n \times n$ matrix where  \\ 
$M_{ij} = \texttt{IsThereBijection}(p_i, S_{p_j})$, where the function \texttt{IsThereBijection}
takes in a prime $p$ and a group $G$ and outputs 1 if there is a bijection between \\
$f:\text{Irr}_{p'}(G) \rightarrow \text{Irr}_{p'}(\textbf{N}_{G}(P))$
such that $f(\chi)(1)$ divides $\chi(1)$ for all $\chi \in \irr G$ and 0 otherwise. Then
notice that there is always a bijection for any entry in the lower triangular side of $M$.

\medskip

For example, for $n = 8$ and $P = \{ 2, 3, 5, 7, 11, 13, 17, 19 \}$, the matrix $M$ would
look like Table~\ref{symmetrictable}. The cells highlighted in green correspond to trivial cases,
as $p$ does not divide the order of the group $G$, and hence the Sylow $p$-subgroup of $G$ is
the identity subgroup. That would mean the normalizer is $G$, and the bijection is simply the identity function. However, the cells highlighted in blue present an interest case which lead
to the following conjecture.

\begin{conjecture}\label{conjsymm}
Let $p$ be prime, and let $S_p$ be the symmetric group on $p$ elements. Then there exists a bijection
$$f:{\rm Irr}_{p'}(S_p) \rightarrow {\rm Irr}_{p'}(\textbf{N}_{S_p}(P)) $$
such that 
$$f(\chi)(1)\, \text{ divides } \, \chi(1) \quad \forall \chi \in \irr{S_p}.$$
\end{conjecture} 

The natural step after seeing the results from Conjecture~\ref{conjsymm} is to generalize
to see if the result holds not only for $S_p$, but see if there is a bijection
for any $S_{p^a}$, where $a \in \mathbb{N}$. The problem with this investigation is computational power. Even
with the optimized algorithm the symmetric group grows too fast for a thorough analysis.
However, for the first few primes and first few powers there is in fact a bijection. 
The author has not found a counterexample yet. In light of this, we pose the following
generalization of Conjecture~\ref{conjsymm}:

\begin{conjecture}\label{conjsymm2}
Let $p$ be prime, let $a \in \mathbb{N}$ and let $S_{p^a}$ be the symmetric group on $p^a$ elements. Then there exists a bijection
$$f:{\rm Irr}_{p'}(S_{p^a}) \rightarrow {\rm Irr}_{p'}(\textbf{N}_{S_{p^a}}(P)) $$
such that 
$$f(\chi)(1)\, \text{ divides } \, \chi(1) \quad \forall \chi \in \irr{S_{p^a}}.$$
\end{conjecture}

As a final remark, let us mention
that we also have checked our conjecture in the so called groups of Lie type that are
stored in \texttt{GAP}, as $\text{GL}_n(q)$, $\text{SL}_n(q)$, etc. (for small values of $n$ and $q$), and it also seems that there will exist bijections with the divisibility property whenever the prime $p$ divides $q$.


\begin{thebibliography}{9}
\bibitem{GAP4}
The GAP~Group, \emph{GAP -- Groups, Algorithms, and Programming, Version 4.8.8}; 2017,
\verb+(https://www.gap-system.org)+
 
\bibitem{Is}
I. M. Isaacs, `{\it Character Theory of Finite Groups}',
AMS-Chelsea, Providence, 2006.
 
\bibitem{McK72}
J. McKay, Irreducible representations of odd degree, 
{\sl J.  Algebra} {\bf 20} (1972), 416--418.
\end{thebibliography}
\end{document}